\documentclass[11pt]{article}
\author{S\'everine Fiedler-Le Touz\'e}
\usepackage{amsfonts, amssymb, epsfig}
\newtheorem{theorem}{Theorem}
\newtheorem{conjecture}{Conjecture}
\newtheorem{definition}{Definition}
\newtheorem{lemma}{Lemma}

\title{M-curves of degree 9 with deep nests}
\begin{document}
\maketitle
\begin{abstract}
The first part of Hilbert's sixteenth problem deals with the classification
of the isotopy types realizable by real plane algebraic curves of given degree $m$.
For $m \geq 8$, one restricts the study to the case of the $M$-curves. For $m=9$, the
classification is still wide open. We say that an $M$-curve of degree 9 has a deep nest if it
has a nest of depth 3. In the present paper, we prohibit 10 isotopy types with deep nest and no outer ovals.
\end{abstract}

\section{Introduction}
Let $A$ be a real algebraic non-singular plane curve of degree $m$.
Its complex part $\mathbb{C}A \subset \mathbb{C}P^2$ is a Riemannian surface
of genus $g=(m-1)(m-2)/2$; its real part $\mathbb{R}A \subset \mathbb{R}P^2$
is a collection of $L \leq g+1$ circles embedded in $\mathbb{R}P^2$.
If $L=g+1$, we say that A is an M-curve.A circle embedded in
$\mathbb{R}P^2$ is called {\em oval\/} or {\em pseudo-line\/} depending
on whether it realizes the class 0 or 1 of $H_1(\mathbb{R}P^2)$.If $m$ is even,
the $L$ components of $\mathbb{R}A$ are ovals; if $m$ is odd,    
$\mathbb{R}A$ contains exactly one pseudo-line,which will be denoted by
$\mathcal{J}$. An oval separates $\mathbb{R}P^2$  into a M\"obius band and
a disc. The latter is called the {\em interior\/} of the oval.
An oval of $\mathbb{R}A$ is {\em empty\/} if its interior contains no other
oval. One calls {\em outer oval\/} an oval that is surrounded by no other
oval. 
Two ovals form an {\em injective pair\/} if one of them lies in the
interior of the other one. We call {\em nest of depth $d+1$\/}
a configuration of ovals $(O_0, O_1, \dots, O_d)$ such that 
$O_i$ lies in the interior of $O_j$ for all pairs $i, j$, with $j > i$.
An $M$-curve of degree $m=2k$ or $m=2k+1$ is said to have a {\em deep nest\/}
if it has a nest of depth $k-1$.Notice that for a curve with deep nest, all 
of the other ovals are empty by Bezout's theorem with an auxiliary line.

Let us call the isotopy type of  $\mathbb{R}A \subset \mathbb{R}P^2$ the
{\em real scheme\/} of $A$; it will be described with the following
notation due to Viro. The symbol $\langle \mathcal{J} \rangle$
stands for a curve consisting in one single pseudo-line; 
$\langle n \rangle$ stands for a curve consisting in $n$ empty ovals.
If $X$ is the symbol for a curve without pseudo-line, $1 \langle X \rangle$
is the curve obtained by adding a new oval, containing all of the
others in its interior. Finally, a curve which is the union of 2 disjoint 
curves $\langle A \rangle$ and $\langle B \rangle$, having the property that
none of the ovals of one curve is contained in an oval of the other curve,
is denoted by $\langle A \amalg B \rangle$.The classification of the real 
schemes which are realizable by $M$-curves of
a given degree in $\mathbb{R}P^2$ is part of Hilbert's sixteenth problem.
This classification is complete up to degree 7 and almost complete in
degree 8. A systematic study of the case $m=9$ has been done, the main 
contribution being due to A. Korchagin. See e.g. \cite{ko1},  \cite{ko3}, 
\cite{ko4}, \cite{ko5}, \cite{or2}
for the constructions, and \cite{ko1}, \cite{ko2}, \cite{fi}, \cite{flt2},
\cite{or3}, \cite{or4} 
for the restrictions.
The main result of the present paper is the prohibition of 10 new schemes.
The proof is an improvement of the classical restriction method: we use Bezout's theorem with auxiliary pencils of rational cubics,
and add to the classical theorems on complex orientations the newer formulas of Orevkov  
\cite{or1} for $M$-curves with deep nests.

Let us briefly recall some facts about complex orientations.
The complex conjugation $conj$ of $\mathbb{C}P^2$ acts on $\mathbb{C}A$
with $\mathbb{R}A$ as fixed points sets. Thus, $\mathbb{C}A \setminus
\mathbb{R}A$ is connected, or splits in 2 homeomorphic halves which are
exchanged by $conj$. In the latter case, we say that $A$ is dividing.
Let us now consider a dividing curve $A$ of degree $m$, and
assume that $\mathbb{C}A$ is oriented canonically.
We choose a half $\mathbb{C}A_+$ of $\mathbb{C}A \setminus \mathbb{R}A$.
The orientation of $\mathbb{C}A_+$ induces an orientation on its
boundary $\mathbb{R}A$. This orientation, which is defined up to complete
reversion, is called {\em complex orientation\/} of $A$.
One can provide all the injective pairs of $\mathbb{R}A$ with a sign as
follows: such a pair is {\em positive\/} if and only if the orientations
of its 2 ovals induce an orientation of the annulus that they bound in
$\mathbb{R}P^2$. Let $\Pi_+$ and $\Pi_-$ be the numbers of positive
and negative injective pairs of $A$. If $A$ has odd degree, each
oval of $\mathbb{R}A$ can be provided with a sign: given an oval $O$
of $\mathbb{R}A$, consider the M\"obius band $\mathcal{M}$
obtained by cutting away the interior of $O$ from $\mathbb{R}P^2$.
The classes $[O]$ and $[2\mathcal{J}]$ of $H_1(\mathcal{M})$ either
coincide or are opposite. In the first case, we say that $O$ is
negative; otherwise $O$ is positive. Let $\Lambda_+$ and $\Lambda_-$
be respectively the numbers of positive and negative ovals of
$\mathbb{R}A$. 
The {\em complex scheme\/} of $A$ is obtained by
enriching the real scheme with the complex orientation:
let e.g. $A$ have real scheme $\langle \mathcal{J} \amalg 1 \langle
\alpha \rangle \amalg \beta \rangle$. The complex scheme of $A$ is
encoded by $\langle \mathcal{J} \amalg 1_{\epsilon} \langle \alpha_+ \amalg
\alpha_- \rangle \amalg \beta_+ \amalg \beta_- \rangle$ where
$\epsilon \in \{ +,- \}$ is the sign of the non-empty oval;
$\alpha_+, \alpha_-$ are the numbers of positive and negative ovals
among the $\alpha$; $\beta_+, \beta_-$ are the numbers of positive
and negative ovals among the $\beta$ (remember that all signs are defined
with respect to the orientation of the pseudo-line $\mathcal{J}$).

\begin{description}
\item[Rokhlin-Mishachev formula:]
{\em If\/} $m = 2k + 1$, {\em then\/}
\begin{displaymath}
2(\Pi_+ - \Pi_-) + (\Lambda_+ - \Lambda_-) = L - 1 - k(k+1)
\end{displaymath}
\item[Fiedler theorem:]
{\em Let $\mathcal{L}_t = \{L_t, t \in [0,1]\}$ be a pencil of real lines
based in a point $P$ of $\mathbb{R}P^2$. Consider two lines $L_{t_1}$
and $L_{t_2}$ of $\mathcal{L}_t$, which are tangent to $\mathbb{R}A$ in
two points $P_1$ and $P_2$, such that $P_1$ and $P_2$ are related by a pair
of conjugated imaginary arcs in $\mathbb{C}A \cap (\bigcup L_t)$.

Orient $L_{t_1}$ coherently to $\mathbb{R}A$ in $P_1$, and transport
this orientation through $\mathcal{L}_t$ to $L_{t_2}$.
Then this orientation of $L_{t_2}$ is compatible to that of $\mathbb{R}A$
in $P_2$.\/}
\item[Orevkov formulas for odd degree]
{\em Let $A$ be an M-curve of degree $m=2k+1$, with a {\em deep nest.\/}
Let $l_+$ and $l_-$ be respectively the numbers of positive and negative non-empty ovals;
$\lambda_+$ and $\lambda_-$ be respectively the numbers of positive and negative
empty ovals.
Let $\pi_s^S, S, s \in \{+, -\}$ be the number of pairs $(O, o)$ where $o$ is an empty oval
surrounded by $O$ and $(S, s)$ are the signs of $(O, o)$.
Then:\/}

$\pi_-^+ - \pi_+^+ = (l_+)^2$, and  $\pi_+^- - \pi_-^- + (\lambda_+ - \lambda_-)/2 = (l_-)^2 + l_-$
 
\end{description}

\section{Restrictions}

\subsection{Definitions and results}
Let $C_9$ be an $M$-curve of degree 9.Given an empty oval $X$ of $C_9$, we 
often will have to consider one point
chosen in the interior of $X$. For simplicity, we shall call this point
also $X$. We denote the pencil of lines based in $X$ by $\mathcal{F}_X$.
Let $X, Y$ be 2 empty ovals of $C_9$. We shall denote by $[XY]$ (resp. $[XY]'$)
the segment of line $XY$ that cuts $\mathcal{J}$ an even (resp. an odd)
number of times. We say that $[XY]$ is the {\em principal segment\/} determined by $X, Y$.
Let $X, Y, Z$ be three empty ovals of $C_9$. Corresponding three
points $X, Y$ and $Z$ determine 4 triangles of $\mathbb{R}P^2$. 
We will call {\em principal triangle\/} and denote by $XYZ$ the triangle
bounded by the segments $[XY]$, $[YZ]$ and $[XZ]$.
We denote respectively by $\mathcal{F}_X: Y \to Z$ and
$\mathcal{F}'_X: Y \to Z$ the pencils of lines based in $X$ that
sweep out the segments $[YZ]$ and $[YZ]'$. We denote by $\mathcal{F}_X^{\epsilon}: Y \to Y$ 
$\epsilon \in \{ \pm \}$ the complete pencils rotating in either direction, $\epsilon = +$ being chosen for
the counter-clockwise direction.
If $Y_1, \dots, Y_n$ are empty ovals met successively by  $\mathcal{F}_X^{\epsilon}: Y_1 \to Y_1$,
we say that this pencil has one {\em $\mathcal{J}$-jump\/}  between 2 consecutive ovals, say $Y_i, Y_{i+1}$
if  $\mathcal{F}_X^{\epsilon}: Y_i \to Y_{i+}$ sweeps out the segment $[Y_iY_{i+1}]'$.
Notice that the total number of $\mathcal{J}$-jumps of a complete pencil over a sequence of ovals is
always odd and does not depend on the choices of $\epsilon$ and $Y_1$.
Denote by $(L_t), t \in [0, 1]$ the lines of the pencil  $\mathcal{F}_X^{\epsilon}: Y_1 \to Y_n$.
We say that $Y_1, \dots, Y_n$ form a {\em Fiedler chain\/} with respect to this pencil if each $Y_i, i = 1, \dots, n-1$
is connected to $Y_{i+1}$ by a pair of imaginary arcs of $\mathbb{C}A \cap (\bigcup L_t)$.
Notice that by Fiedler's theorem, the signs of $Y_i$ and $Y_{i+1}$ coincide if and only if the pencil 
$\mathcal{F}_X^{\epsilon}: Y_i \to Y_{i+1}$ has a $\mathcal{J}$-jump.  

An ordered group of empty ovals $A_1, \ldots, A_n$ of $C_9$
{\em lies in a convex position\/} if for each triple
$A_i, A_j, A_k$, the principal triangle $A_iA_jA_k$ does not contain
any other oval of the group and $A_1, \ldots, A_n$ are the successive
vertices of $\bigcup A_iA_jA_k$ (the {\em convex hull of the group\/}).

Let $C_9$ be an $M$-curve with real scheme  
$\langle \mathcal{J} \amalg \alpha \amalg 1 \langle \beta \amalg 1 \langle \gamma \rangle \rangle \rangle$.
We shall call the ovals of the groups $\alpha$, $\beta$ and $\gamma$ respectively
{\em outer\/}, {\em median\/} and {\em inner\/} ovals.
Notice that the inner ovals must lie in convex position. Otherwise, let $A, B, C, D$ be 4
inner ovals such that $D$ lies inside of the principal triangle $ABC$.Then, for any choice
of a fifth empty oval $E$, the conic through $A, B, C, D, E$ cuts $C_9$ in more than 18 points.
This is a contradiction. 

\begin{definition}
We say that $C_9$ has a {\em $O_1$-jump\/} if there exist 2 median ovals
$A, C$ and 2 inner ovals $B, D$, such that the line $AC$ separates $B$ from $D$ in
$Int(O_1)$. 
\end{definition}

Notice, again with auxiliary conics, that if $C_9$ has a $O_1$-jump determined by $A, B, C, D$, then
the principal segment $[AC]$ cuts $O_1$, as shown in Figure 1.
The lines $(AB), (BC), (CD), (DA)$ give rise to 3 quadrangles and 4
triangles. All of the remaining empty ovals lie in the union of the 4 triangles.

\begin{figure}
\centering \psfig{file=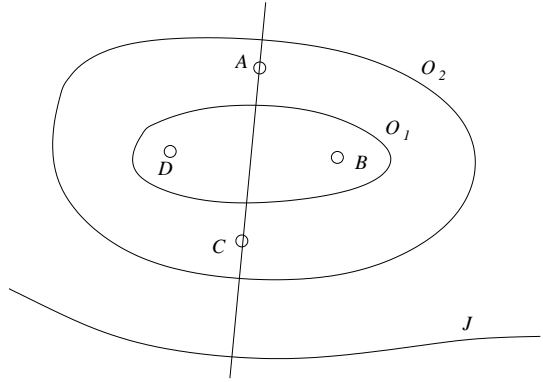}
\caption{$O_1$-jump}
\end{figure}

\begin{lemma}
If $C_9$ has a median oval $A_1$ such that the pencil $\mathcal{F}_{A_1}$ sweeping out $O_1$
meets successively empty ovals $A_2, A_3, \dots, A_{2J+2}$, where the $A_i$, $i$ even (resp. odd)
are inner (resp. median) ovals.
Then:

(1) The $A_i$ $i=1, \dots, 2J+2$ lie in convex position.  

(2) Let $A_i, A_j, A_k, A_l$ be 4 consecutive ovals in the cyclic ordering, we denote by $T_j$
the triangle determined by the lines $(A_iA_j)(A_jA_k)(A_kA_l)$ with edge $[A_jA_k]$, that does not
intersect the convex hull of $A_1, A_2, \dots, A_{2J+2}$.
All of the remaining empty ovals lie in $\bigcup T_j$, $j= 1, \dots, 2J+2$; if an oval $A$ lies in $T_j \cap Int(O_2)$, $A$ is not
separated in $T_j$ from the edge $[A_jA_k]$ by $\mathcal{J}$. 

(3) There is a natural cyclic ordering of the empty ovals. This ordering is given by:
the complete pencils of lines $\mathcal{F}_B$, for any inner oval $B$, and
the pencils $\mathcal{F}_A$, where $A$ is a median oval in some triangle $T$, sweeping out the $2J$ triangles $T_j$
having no common median vertex with $T$.
The ovals in the triangles $T_1, T_2, \dots, T_{2J+2}$ appear successively in the cyclic ordering.

(4) If $\alpha = 0$, then $\lambda_+ - \lambda_- = 0$
\end{lemma}

{\em Proof\/} The points (1), (2), (3) are easily proven using Bezout's theorem with conics. 
Notice that if $\alpha = 0$, the pencils of lines $\mathcal{F}_{A_i}: A_j \to T_j \to A_k$ have no $\mathcal{J}$-jumps, 
and give rise to a closed Fiedler chain involving all of the empty ovals, hence (4) follows. $\Box$

\begin{figure}
\centering \psfig{file=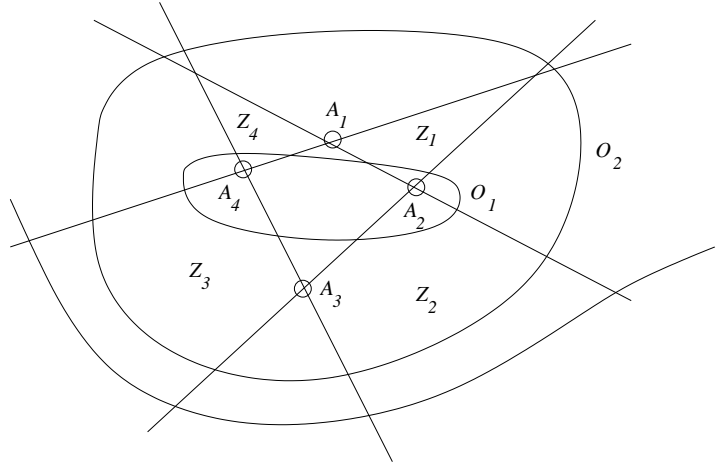}
\caption{$C_9$ with $J=1$}
\end{figure}

\begin{definition}
Let $C_9$ be an $M$-curve with deep nest, and no outer ovals.
The curve $C_9$ has $J$ {\em $O_1$-jumps\/} with distribution $(l_1, \dots, l_{2J+2})$ if $C_9$ verifies the
condition of Lemma 1 and $l_i, i=1, \dots, 2J+2$ are the cardinals of the successive groups of inner
and median ovals in the cyclic ordering.
\end{definition}

The case $J=1$ is illustrated in Figure 2.
Notice that the number of $O_1$-jumps and their distribution is a rigid isotopy invariant of $C_9$.

\begin{theorem}
The real schemes 
$\langle \mathcal{J} \amalg 1 \langle \beta \amalg 1 \langle \gamma \rangle \rangle \rangle$,
with $\beta + \gamma = 26$ and $\beta, \gamma$ odd are not realizable by $M$-curves of degree 9.
\end{theorem}

Among the 13 admissible schemes of that form, the 3 with $\beta = 1, 3$ and $25$ had already been forbidden in \cite{ko2}.

\begin{theorem}
Let $C_9$ be an $M$-curve with real scheme
$\langle \mathcal{J} \amalg 1 \langle \beta \amalg 1 \langle \gamma \rangle \rangle \rangle$, with even $\beta, \gamma$.
Then:
\begin{enumerate}
\item
if $C_9$ has $O_1$-jumps, then the complex scheme of $C_9$ is

$\langle \mathcal{J} \amalg 1_- \langle (\frac{\beta - 6}{2})_+ \amalg (\frac{\beta + 6}{2})_-
\amalg 1_- \langle (\frac{\gamma + 6}{2})_+ \amalg (\frac{\gamma - 6}{2})_-
\rangle \rangle \rangle$ or

$\langle \mathcal{J} \amalg 1_+ \langle (\frac{\beta + 4}{2})_+ \amalg (\frac{\beta - 4}{2})_-
\amalg 1_+ \langle (\frac{\gamma - 4}{2})_+ \amalg (\frac{\gamma + 4}{2})_-
\rangle \rangle \rangle$.

\item
otherwise, the complex scheme of $C_9$ is

$\langle \mathcal{J} \amalg 1_- \langle (\frac{\beta + 4}{2})_+ \amalg (\frac{\beta - 4}{2})_-
\amalg 1_- \langle (\frac{\gamma}{2})_+ \amalg (\frac{\gamma}{2})_-
\rangle \rangle \rangle$.
\end{enumerate}
\end{theorem}

Notice that by Theorem 2, the real scheme $\langle \mathcal{J} \amalg 1 \langle 2 \amalg
1 \langle 24 \rangle \rangle \rangle$ is not realizable. This real scheme was already forbidden in \cite{ko2}.
S. Orevkov constructed $M$-curves realizing the 10 real schemes 
$\langle \mathcal{J} \amalg 1 \langle \beta \amalg 1 \langle \gamma \rangle \rangle \rangle$, with $\beta, \gamma$ even,
with $4 \leq \beta \leq 22$. Each of these curves has 3 $O_1$-jumps, they realize the complex schemes 
$\langle \mathcal{J} \amalg 1_+ \langle (\frac{\beta + 4}{2})_+ \amalg (\frac{\beta - 4}{2})_-
\amalg 1_+ \langle (\frac{\gamma - 4}{2})_+ \amalg (\frac{\gamma + 4}{2})_-
\rangle \rangle \rangle$ (private communication).

\subsection{Lemmas}

Assume there exists an $M$-curve $C_9$ with real scheme
$\langle \mathcal{J} \amalg 1 \langle \beta \amalg
1 \langle \gamma \rangle \rangle \rangle$.

We repeat hereafter the arguments from  \cite{ko2} for the case
$\beta = 0$: let $A, B, C$ be 3 inner ovals. The pencils of lines $\mathcal{F}_A: B \to C$,
$\mathcal{F}_B: C \to A$, and $\mathcal{F}_C: A \to B$ have no $\mathcal{J}$-jumps and give rise to a closed 
Fiedler chain involving all of the inner ovals. Therefore, $\lambda_+ - \lambda_- = 0$, $\Pi_+ - \Pi_- = \pm 1$ and
$\Lambda_+ - \Lambda_- \in \{ 0, 2, -2 \}$. This contradicts the Rokhlin-Mishachev formula.
The real scheme $\langle \mathcal{J} \amalg 1 \langle 1 \langle 26 \rangle \rangle \rangle$ is not realizable.

Let now $\beta > 0$. Assume $C_9$ has $O_1$-jumps. Let $\epsilon_1, \epsilon_2$ and 
$n\epsilon_3$ be respectively the contributions of $O_1$, $O_2$ and of the inner ovals to
$\Lambda_+ - \Lambda_-$ (where $\epsilon_1, \epsilon_2, \epsilon_3 \in \{+1, -1\}$, and $n$ is non-negative). 
By Lemma 1 (4) and the Rokhlin-Mishachev formula, one must have:
$2(-\epsilon_1\epsilon_2 - n\epsilon_1\epsilon_3) + \epsilon_1 + \epsilon_2 = 8$.
There are 4 possible solutions:

$\epsilon_1 = \epsilon_2 = -1$, $\epsilon_3 = 1$, $n=6$;

$\epsilon_1 = -1$, $\epsilon_2 = \epsilon_3 = 1$, $n=3$;

$\epsilon_1 = 1$, $\epsilon_2 = \epsilon_3 = -1$, $n=3$;

$\epsilon_1 = \epsilon_2 = 1$, $\epsilon_3 = -1$, $n=4$.

In the second case, one has $l_+ = 1$, $\pi_-^+ - \pi_+^+ = 0$;
in the third case, one has  $l_+ = 1$, $\pi_-^+ - \pi_+^+ = 3$. 
This contradicts the first Orevkov formula.
Either of the other 2 cases verifies both Orevkov formulas.
The numbers $n$, $\beta$, $\gamma$ have the same parity. If $C_9$ has $O_1$-jumps, these numbers are even,
and Theorem 2 (1) is proven.

\begin{lemma} 
Let $C_9$ be an $M$-curve with real scheme $\langle \mathcal{J} \amalg 1 \langle \beta \amalg
1 \langle \gamma \rangle \rangle \rangle$, without $O_1$-jumps and such that $\gamma \geq 2$ and 
$\beta \geq 1$.
Let $A$ be a median oval, and $B, C$ be the extreme inner ovals met by the pencil $\mathcal{F}_A$
sweeping out $O_1$. 
The complete pencil $\mathcal{F}_C$ gives rise to a cyclic Fiedler chain 
involving all other empty ovals. In the corresponding cyclic ordering, all of the inner ovals are consecutive.
\end{lemma}

{\em Proof:\/}
Let $T_1, \dots, T_4$ be the 4 triangles $ABC$, $T_1$ and $T_2$ being the 2 ones
with edge $[BC]$. The inner ovals lie in $Z'_1 \cup Z'_2$, where $Z'_i = T_i \cap Int(O_1)$, $i=1, 2$.
One of the zones  $Z'_i, i=1, 2$ is empty. Indeed, assume either zone contains an oval $E_i$.
Then the conic through $A, B, C, E_1, E_2$ cuts $C_9$ in more than 18 points. Contradiction.
Consider the pencil of lines $\mathcal{F}_C$ starting at $B$ and sweeping out the non-empty zone $Z'_i, i \in \{1, 2\}$.
Let $D$ be the last inner oval met by this pencil. The pencil $\mathcal{F}_C: B \to D$ meets no median oval.
Indeed, assume that a median oval $E$ is met by this pencil. Then, the conic $EBADC$ cuts $C_9$
in more than 18 points. Contradiction. The complete pencil $\mathcal{F}_C$ gives rise to a cyclic 
Fiedler chain involving all other empty ovals. This chain splits into 2 consecutive subchains formed respectively by the
inner and the median ovals. $\Box$

\begin{lemma}
Let $C_9$ be an $M$-curve of degree 9 with deep nest, and $1$ be any one of the inner ovals.
Assume there exist five other empty ovals $2, 3, 4, 5, 6$ met successively by the pencil of lines 
$\mathcal{F}_1^+$, and such that there is a $\mathcal{J}$-jump between any two successive ovals in the cyclic
ordering. Denote by $\mathcal{F}_{1123456}$ the pencil of rational cubics through $1, \dots, 6$ with double point at $1$.
Then, up to cyclic permutation of $2, \dots, 6$, the sequence of singular (i.e. reducible) cubics of $\mathcal{F}_{1123456}$ is:

$16 \cup 14523, 14 \cup 12356, 12 \cup 14365, 15 \cup 12643, 13 \cup 15426$ or

$12 \cup 14365, 15 \cup 12643, 16 \cup 12543, 13 \cup 12456, 14 \cup 12356$ or

$16 \cup 15234, 14 \cup 15326, 15 \cup 13264, 13 \cup 14265, 12 \cup 14365$
\end{lemma}

See Figures 9, 10, 11.

{\em Proof:\/}
Notice first that if $\mathcal{F}_1^\epsilon: X, Y, Z, X$ has a $\mathcal{J}$-jump between $X$ and $Y$ and between
$Y$ and $Z$, then: $\mathcal{F}_1^\epsilon$ has also a $\mathcal{J}$-jump between $Z$ and $X$, and $1$ lies in the
principal triangle $XYZ$. Thus, $1 \in 234 \cap 345 \cap 456 \cap 562 \cap 623$. 
As no three points among $1, \dots, 6$ can be on a line, these five principal triangles must have a 2-dimensional intersection. 

{\em Case 1:\/} The five points lie in convex position. 
Consider a pair of points, that are consecutive for $\mathcal{F}_1^+$, say $2, 3$, and assume these points also consecutive for 
the convex cyclic ordering. This cyclic ordering is $2, 3, X, Y, Z$, with $X, Y, Z \in \{4, 5, 6\}$.
If $X, Y, Z = 6, 5, 4$ or $5, 6, 4$, then $234 \cap 345 = [34]$; 
if $X, Y, Z = 6, 4, 5$ or $5, 4, 6$, then $632 \cap 456 = 6$;
if $X, Y, Z = 4, 6, 5$ or $4, 5, 6$, then $234 \cap 456 = 4$. Contradiction.
Then, the only possible convex cyclic ordering of the five points is $2, 4, 6, 3, 5$.
Choose a line at infinity $L$ that does not cut the convex hull of the points. In the oriented affine plane
$\mathbb{R}P^2 \setminus L$, there are a priori two possibilities for the {\em positive\/} cyclic convex ordering of the five points:
$2, 4, 6, 3, 5$ and $2, 5, 3, 6, 4$.
As  $\mathcal{F}_1^+$ sweeps out successively $2, 3, 4, 5, 6$, the first of these possibilities is realized.
(see Figure 3).

{\em Case 2\/:} One of the points, say $2$, lie in the convex hull of the other four. 
Assume $3$ and $4$ are consecutive in the cyclic convex ordering. Then, $234 \cap 562 = 2$. Contradiction.
There are a priori two possibilities for the positive cyclic convex ordering of the four points: $3, 5, 4, 6$ or $3, 6, 4, 5$.
In the first case, the quadrangle $3546$ is divided into four triangles: $T_1 = 345 \cap 356$, $T_2 = 356 \cap 346$, 
$T_3 = 346 \cap 456$ and $T_4 = 345 \cap 456$. If $2 \in T_1 \cup T_2$, then $623 \cap 456 = 6$; if $2 \in T_3$, then 
$623 \cap 345 = 3$. Contradiction. Hence, $2 \in T_4$. 
In the second case, the quadrangle $3645$ is divided into four triangles: $T_1 = 346 \cap 356$, $T_2 = 356 \cap 345$,
$T_3 = 345 \cap 456$, $T_4 = 346 \cap 456$. If $2 \in T_1 \cup T_4$, then $632 \cap 345 = 3$; if $2 \in T_2$, then
$623 \cap 456 = 6$. If $2 \in T_3$, then $1 \in 234 \cap 562$ and the pencil  
$\mathcal{F}_1^+$ sweeps out successively $6, 5, 4, 3, 2$.
This is a contradiction. Thus, $2 \in T_2$ (see Figure 4).
 
{\em Case 3\/}: Two of the points lie in the principal triangle determined by the other three. 
Assume the two points are not consecutive for $\mathcal{F}_1^+$, say these points are $4, 6$. 
Then $234 \cap 345 = [34]$. Contradiction. Thus, up to cyclic permutation of $2, 3, \dots, 6$, one can choose 
$4, 5$ as interior points.
There are a priori two possibilities for the positive convex ordering of the other three points: $2, 6, 3$ and $2, 3, 6$.
As $1 \in 632$ and $\mathcal{F}_1^+$ sweeps out successively $6, 2, 3$, the latter possibility is excluded.
The triangle $632$ is divided in six triangles $T_i, i=1, \dots, 6$ by the lines $42, 43$ and $46$, such that:
$T_1 \cup T_2 = 346$, $T_3 \cup T_4 = 234$,  $T_5 \cup T_6 = 246$; and 
$T_1, T_6$ have $4, 6$ as common vertices, $T_2, T_3$ have $4, 3$ as common vertices,
$T_4, T_5$ have $4, 2$ as common vertices. If $5 \in T_1 \cup T_2 \cup T_4 \cup T_5$, then $345 \cap 456 = [45]$;
if $5 \in T_6$, then $456 \cap 562 = [56]$. Contradiction. 
One has $5 \in T_3$ (see Figure 5).
 
Perform a Cremona transformation $cr: (x_0; x_1; x_2) \to (x_1x_2; x_0x_2; x_0x_1)$ with base points $1, 5, 4$. 
Let us denote the respective images of the lines $14, 15, 45$ by $5, 4, 1$. For the other points, we shall use
the same notation as before $cr$. After $cr$, consider the pencil of conics $\mathcal{F}_{1236}$. 
This pencil has 3 singular conics: the double lines $12 \cup 36$, $13 \cup 26$ and $16 \cup 23$, they are shown in
Figures 6, 7 and 8. The position of the pencil $\mathcal{F}_{1236}$ with respect to the base lines is
uniquely determined, it is shown in the left-hand side of Figures 9, 10, and 11. In cases 1, 2, this is obvious; in case 3, 
both points $4, 5$ are swept out in the portion $16 \cup 23 \to 13 \cup 26$ and we have to find out which one is met first.
The conic $16235$ of $\mathcal{F}_{1236}$ is the image of a conic $15326$ (see Figure 5), thus $4$ lies outside
of $16235$. So, $5$ is swept out before $4$ by $\mathcal{F}_{1236}$.
Perform the Cremona transformation back. The pencil of conics $\mathcal{F}_{1236}$ is mapped onto the pencil of rational cubics 
$\mathcal{F}_{1123456}$ (right-hand side of Figures 9, 10, 11). $\Box$

\begin{lemma}
Let $C_9$ be an $M$-curve of degree 9 with deep nest, and denote by $1$ any one of the inner ovals.
The complete pencil of lines $\mathcal{F}_1$ has at most $3$ $\mathcal{J}$-jumps over the sequence of median ovals.
\end{lemma}

{\em Proof\/}
Assume there exists an inner oval $1$ of $C_9$ such that $\mathcal{F}_1$ has $5$ $\mathcal{J}$-jumps
over the median ovals. Let $2, \dots, 6$ be median ovals such that $1, 2, \dots, 6$ verify the conditions of Lemma 3.
Notice that in the affine plane of Figures 9, 10, 11, a segment of cubic 
connecting two points among $1, \dots, 6$ goes over infinity if and only if it cuts $\mathcal{J}$.
In each of the 3 cases, the cubics of the pencil $\mathcal{F}_{1123456}$
intersect all $O_1 \cup O_2 \cup \{1, 2, 3, 4, 5, 6\} \cup \mathcal{J}$ at 27 points, so there is no
possibility to sweep out the other empty ovals of $C_9$. Contradiction. $\Box$

\begin{figure}
\centering \psfig{file=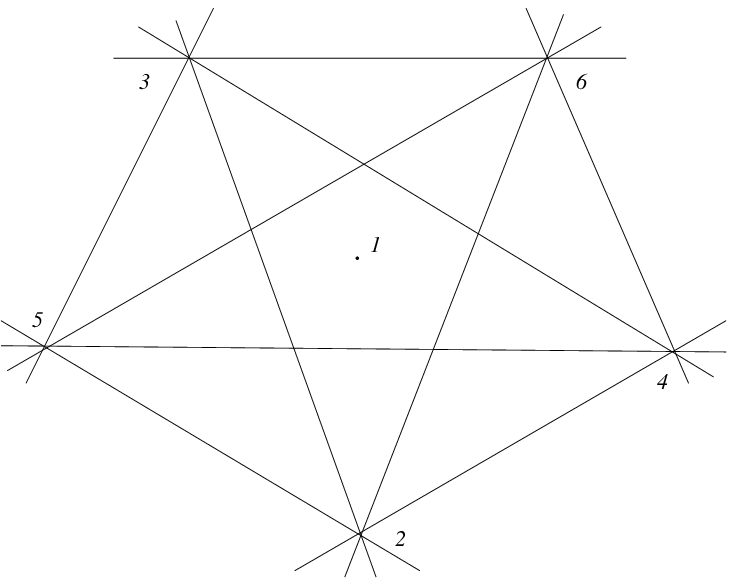}
\caption{Case 1}
\end{figure}

\begin{figure}
\centering \psfig{file=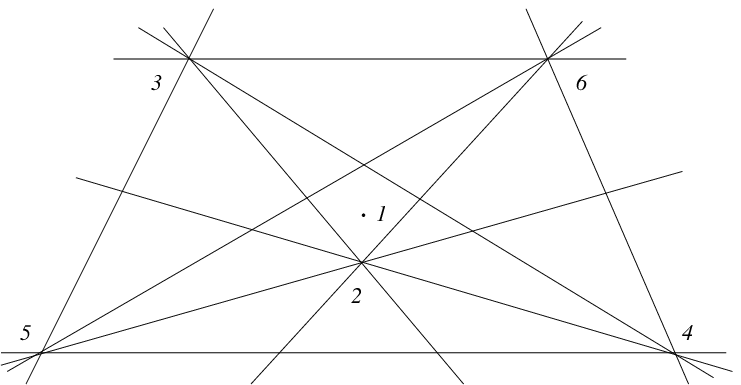}
\caption{Case 2}
\end{figure}

\begin{figure}
\centering \psfig{file=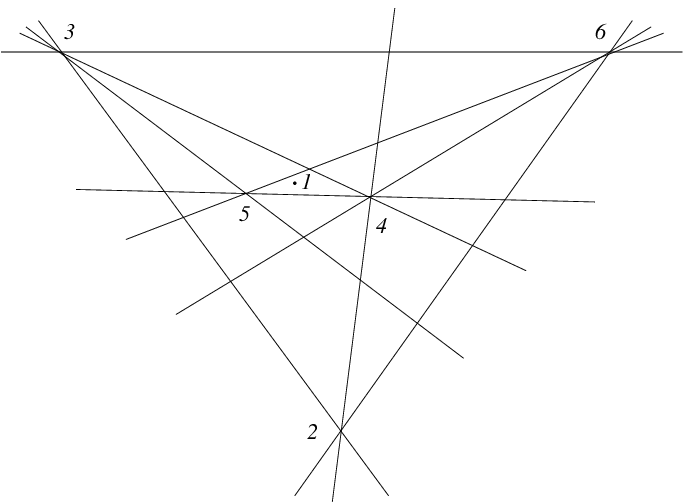}
\caption{Case 3}
\end{figure}

\begin{figure}
\centering \psfig{file=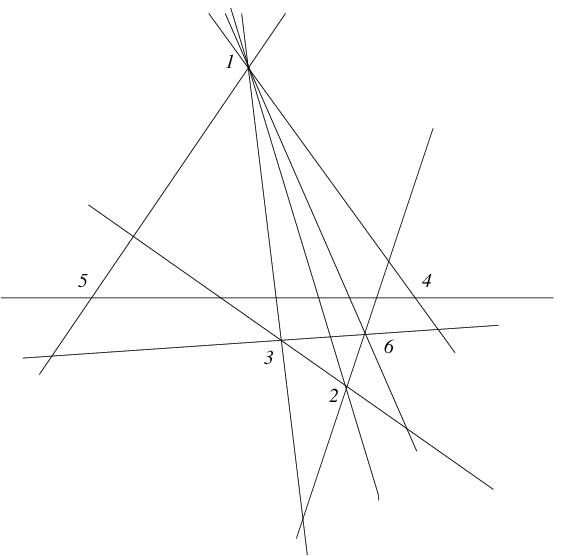}
\caption{The 3 double lines of $\mathcal{F}_{1236}$, case 1}
\end{figure}

\begin{figure}
\centering \psfig{file=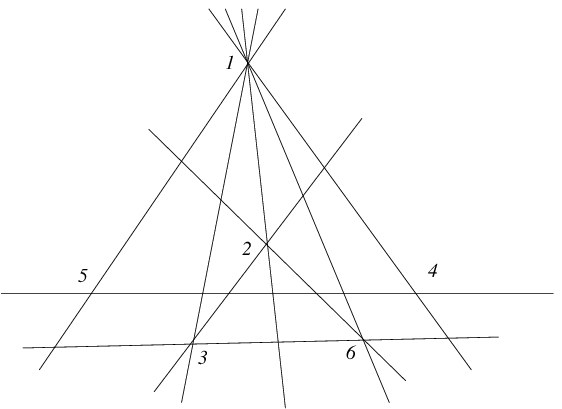}
\caption{The 3 double lines of $\mathcal{F}_{1236}$, case 2}
\end{figure}

\begin{figure}
\centering \psfig{file=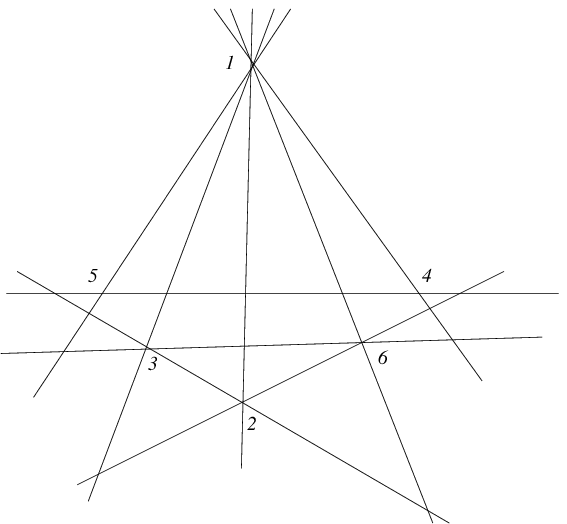}
\caption{The 3 double lines of $\mathcal{F}_{1236}$, case 3}
\end{figure}

\begin{figure}
\centering \psfig{file=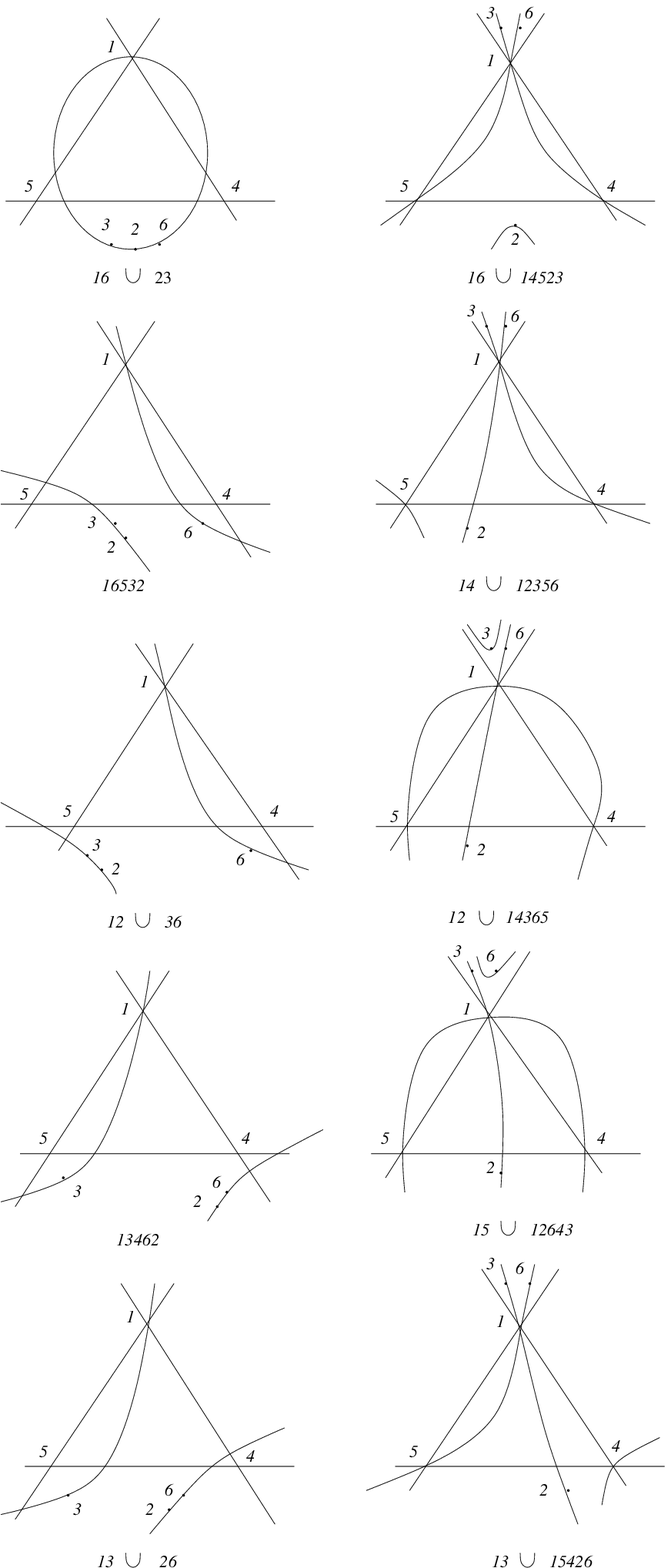}
\caption{$\mathcal{F}_{1123456} = cr^{-1}(\mathcal{F}_{1236})$, case 1}
\end{figure}

\begin{figure}
\centering \psfig{file=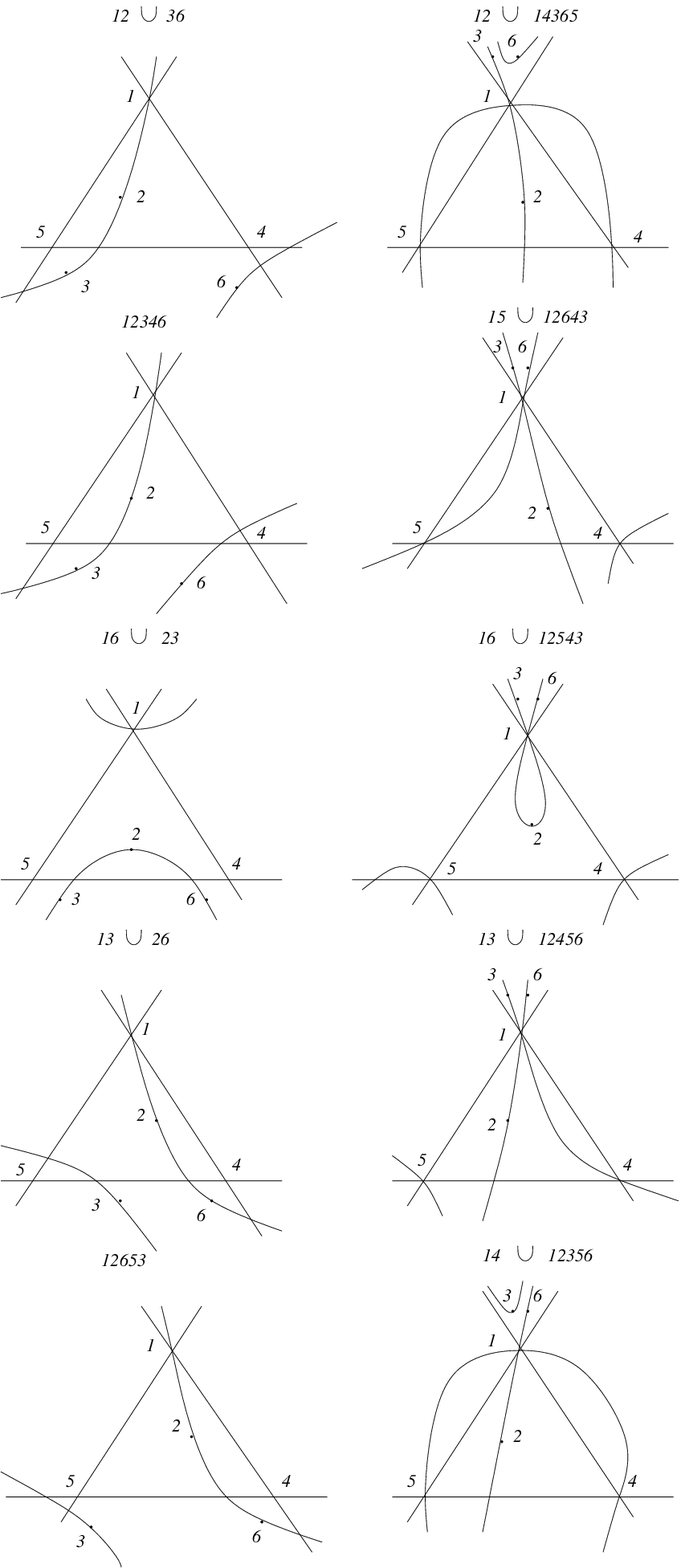}
\caption{$\mathcal{F}_{1123456} = cr^{-1}(\mathcal{F}_{1236})$, case 2}
\end{figure}

\begin{figure}
\centering \psfig{file=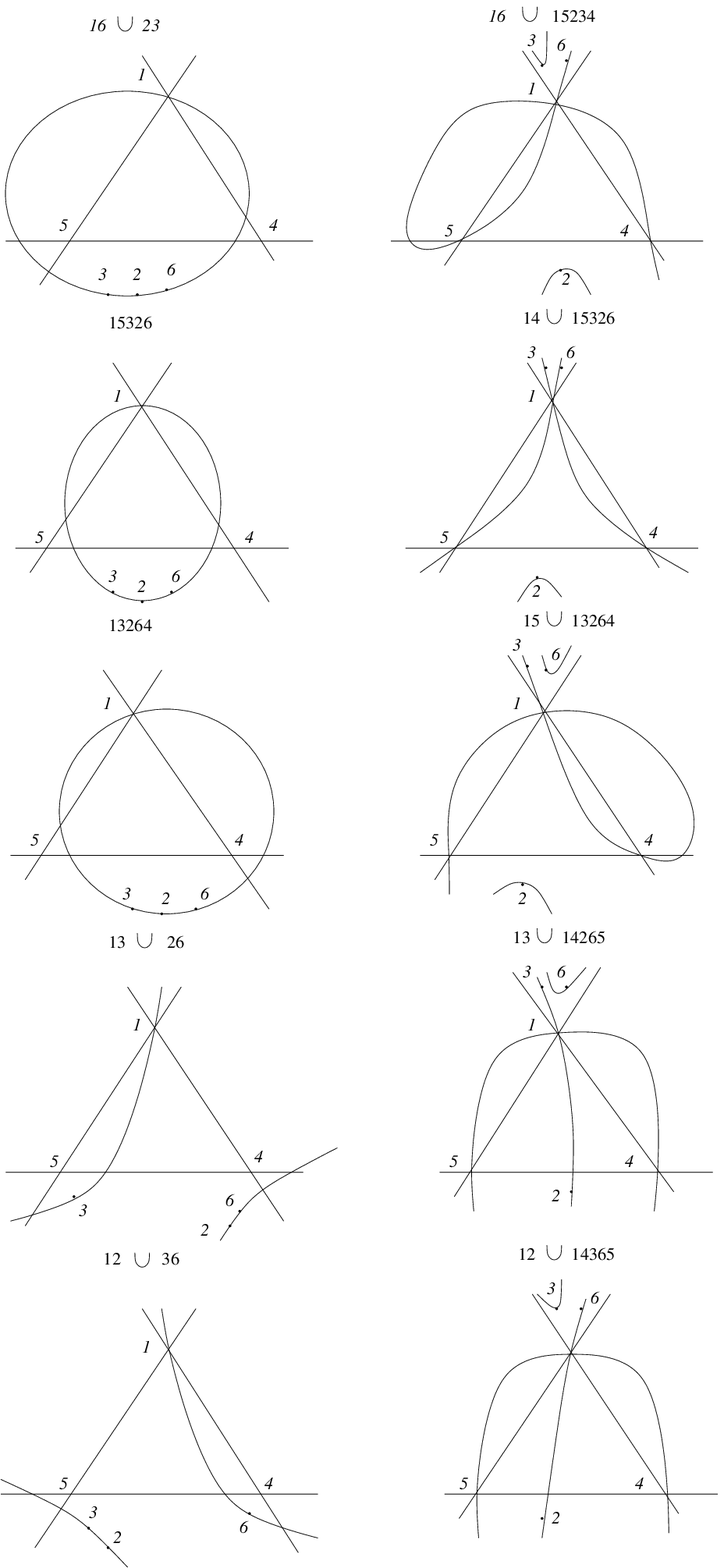}
\caption{$\mathcal{F}_{1123456} = cr^{-1}(\mathcal{F}_{1236})$, case 3}
\end{figure}

\subsection{Proof of Theorems 1 and 2}

Let $C_9$ be an $M$-curve of degree 9 with real scheme  $\langle \mathcal{J} \amalg 1 \langle \beta \amalg
1 \langle \gamma \rangle \rangle \rangle$ and without $O_1$-jumps.
Let $\epsilon_1, \epsilon_2$ and $n\epsilon_3$ be respectively the contributions of $O_1$, $O_2$ and of 
the median ovals to $\Lambda_+ - \Lambda_-$; if $\gamma$ is odd, let $\epsilon_4$ be the contribution of the inner ovals to
$\Lambda_+ - \Lambda_-$ 
(where $\epsilon_1, \epsilon_2, \epsilon_3, \epsilon_4 \in \{+1, -1\}$ and $n$ is non-negative).
Let $A$ be a median oval and $C$ be one of the extreme inner ovals met by $\mathcal{F}_A$.
Applying Lemmas 2 and 4 with the complete pencil of lines $\mathcal{F}_C$, we prove that if $\gamma$ is odd, then $n \in \{1, 3\}$ and
if $\gamma$ is even, $n \in \{0, 2, 4 \}$. 
The Rokhlin-Mishachev formula yields respectively for odd and for even $\gamma$:

\begin{displaymath}
2(-\epsilon_1\epsilon_2 - \epsilon_4\epsilon_2 - \epsilon_4\epsilon_1 - n\epsilon_3\epsilon_2) +
\epsilon_1 + \epsilon_2 + n\epsilon_3 + \epsilon_4 = 8 
\end{displaymath}

\begin{displaymath}
2(-\epsilon_1\epsilon_2 - n\epsilon_3\epsilon_2) + \epsilon_1  + \epsilon_2 + n\epsilon_3 = 8 
\end{displaymath}

For $\gamma$ odd, there is no solution. For $\gamma$ even, there are 2 solutions:

$\epsilon_1 = \epsilon_2 = -1$, $\epsilon_3 = 1$, $n = 4$;

$\epsilon_1 = \epsilon_3 = 1$, $\epsilon_2 = -1$, $n = 2$.

In the latter case, one has: $l_+ = 1$, $\pi_-^+ - \pi_+^+ = 0$. Contradiction with the first Orevkov formula.
In the first case, both Orevkov formulas are verified.
The numbers $n$, $\beta$, $\gamma$ have the same parity, if $C_9$ has no $O_1$-jump, these numbers are even. 
This finishes the proof of Theorem 1 and of Theorem 2 (2). $\Box$

\end{document}